\documentclass[12pt]{amsart}
\usepackage{a4wide, amssymb, graphicx, times}
\usepackage[colorlinks, citecolor=black, filecolor=black, linkcolor=black,
urlcolor=black]{hyperref}

\newtheorem{theorem}{Theorem}[section]
\newtheorem{definition}[theorem]{Definition}
\newtheorem*{theorem*}{Theorem}

\newtheorem{corollary}[theorem]{Corollary}

\newcommand{\la}{\langle}
\newcommand{\ra}{\rangle}
\newcommand{\Z}{\mathcal{Z}}
\newcommand{\C}{\mathbb{C}}

\renewcommand{\S}{\mathbb{S}}

\begin{document}

\title{Equidistribution of the Fekete points on the sphere}
\author{Jordi Marzo}
\author{Joaquim Ortega-Cerd\`a}
\date{\today}
\begin{abstract} 
The Fekete points are the points that maximize a Vandermonde-type determinant
that appears in the polynomial Lagrange interpolation formula. They are well
suited points for interpolation formulas and numerical integration. We prove the
asymptotic equidistribution of the Fekete points in the sphere. The way we
proceed is by showing their connection with other array of points, the
Marcinkiewicz-Zygmund arrays and the interpolating arrays, that have been
studied recently. 
\end{abstract}
\maketitle

\section{Introduction}

For any integer $\ell\ge 0,$ let $\mathcal{H}_{\ell}$ be the space of spherical
harmonics of degree $\ell$ in $\S^{d}$. For any integer $L\ge 0$ we denote the
space of spherical harmonics of degree not exceeding $L$ by    $\Pi_{L}$ These
vector spaces have dimensions 
\[
\dim \Pi_{L}=\frac{d+2L}{d}\binom{d+L-1}{L}=\pi_{L}\simeq L^d.
\]

Let $\{ Q_1^L ,\ldots , Q_{\pi_L}^L \}$  be any basis in $\Pi_{L}$. The
points $\Z_L=\{ z_{L1} ,\ldots , z_{L,\pi_L}  \}$ maximizing the determinant
\[
|\Delta( x_{1} ,\ldots , x_{\pi_L} )|=|\det (Q^L_i(x_{j}))_{i,j}|
\]
are called the Fekete points of degree $L$ for $\S^d$ (this points are sometimes
called extremal fundamental systems of points as in \cite{SloWom04}). They are
not to be confused with the \emph{elliptic} Fekete points which are a system of
points that minimize the potential energy. The extremal fundamental system of 
points are better suited nodes for cubature formulas and for polynomial
interpolation, see \cite{SloWom04} and the references therein.


The geometric properties of the distribution of the Fekete points on the sphere
has been the subject of research, see for instance \cite{Reimer90},
\cite{BosLevWal08} or \cite{SloWom04}. One natural problem is the limiting
distribution of the points as $L\to\infty$. If we denote by $\mu_L= \frac
1{\pi_L} \sum_j\delta_{x_{j}}$ it has been known long ago that for the
\emph{elliptic} Fekete points $\mu_L$ converges vaguely to the uniform
distribution of the sphere even for a wide class of potentials, see
\cite{HarSaf04} for a very nice survey. For the Fekete points in  compacts
$K\subset \C$ this is a classical result. Much less is known in higher
dimensions. On the paper \cite{BerBou08} the authors have found the limiting
distribution in the context of line bundles over complex manifolds. The
techniques there are very different from ours, they rely on a careful study of
the weighted transfinite diameter and its differentiability. We will rather
emphasize the connection of the Fekete points with the Marcinkiewicz-Zygmund
arrays and the interpolating arrays (see below for the definitions). As long as
the density of these arrys is understood we can obtain the equidistribution of
the Fekete points. This is the case of the sphere where we can build on the work
\cite{Marzo07}, where the M-Z arrays and interpolating arrays are studied. The
same approach is being pursued by R.~Berman in line bundles over complex
manifolds.

\subsection{Marcinkiewicz-Zygmund inequalities and interpolation}

We consider arrays of points on the sphere $\S^{d}$ that determine the norm of
the polynomials, and also arrays of points where we are free to interpolate
arbitrary values by polynomials,i.e. interpolating arrays. More precisely, For
any degree $L$ we take $m_{L}$ points in $\S^{d}$
\[
\Z(L)=\{ z_{L j}\in \S^{d}: 1\le j\le m_{L}\}, \quad L\ge 0,
\]
and assume that $m_{L}\to \infty$ as $L\to \infty$. This yields a
triangular array of points $\Z=\{ \Z(L) \}_{L\ge 0}$ in $\S^{d}$.

\begin{definition}\label{def-MZ}
Let $\Z=\{ \Z(L) \}_{L\ge 0}$ be a triangular array with $m_{L}\ge \pi_{L}$ for
all $L$. We call $\Z$ an $L^{p}$-Marcinkiewicz-Zygmund array, denoted by
$L^{p}$-MZ, if there exists a constant $C_{p}>0$ such that for all $L\ge 0$ and
$Q\in \Pi_{L}$,
\begin{equation}\label{ineq-def-MZ}
\frac{C_{p}^{-1}}{\pi_{L}}\sum_{j=1}^{m_{L}}|Q(z_{L j})|^{p} \le
\int_{\S^{d}}|Q(\omega)|^{p}d\sigma(\omega) \le
\frac{C_{p}}{\pi_{L}}\sum_{j=1}^{m_{L}}|Q(z_{Lj})|^{p}, 
\end{equation}
if $1\le p<\infty,$ and 
\[
\sup_{\omega\in \S^{d}}|Q(\omega)|\le C\sup_{j=1,\ldots , m_{L}}|Q(z_{Lj})|,
\] 
when $p=\infty$.
\end{definition}

Then the $L^{p}-$norm in $\S^{d}$ of a polynomial of degree $L$ is comparable to
the discrete version given by the weighted $\ell^{p}-$norm of its restriction to
$\Z(L)$.

\begin{definition}\label{def-interp}
Let $\Z=\{ \Z(L) \}_{L\ge 0}$ be a triangular array with $m_{L}\le \pi_{L}$ for
all $L$. We say that $\Z$ is $L^{p}-$interpolating, if for all arrays $\{
c_{Lj}\}_{L\ge 0, 1\le j\le m_{L}}$ of values such that
\[
\sup_{L\ge 0}\frac{1}{\pi_{L}}\sum_{j=1}^{m_{L}}|c_{Lj}|^{p}<\infty,
\]
there exists a sequence of polynomials $Q_{L}\in \Pi_{L}$ uniformly bounded in
$L^{p}$ such that $Q_{L}(z_{Lj})=c_{Lj},$ $1\le j\le m_{L}$.
\end{definition}

Roughly speaking  in order to recover the $L^{p}-$norm of a polynomial of degree
$L$ from the evaluation at the points in $\Z(L)$ we need a sufficiently big
number of points in $\Z(L)$. Thus, intuitively, the M-Z arrays must have high
density. On the other hand, in an interpolating array it is possible to have a
spherical harmonic of degree at most $L$ attaining some prescribed values on
$\Z(L)$. Intuitively this is possible only when $\Z(L)$ is sparse. 

In dimension one the roots of unity are simultaneously an interpolating and an
M-Z array when $1<p<\infty$. On higher dimension the situation is more delicate.
It has been proved, \cite[Theorem 1.7]{Marzo07} that there are not arrays which
are simultaneously $L^p$-MZ and interpolating when $d>2$ and $p\ne 2$ and most
likely even when $p=2$. We will prove that the Fekete points are a very
reasonable substitute. If we perturb them slightly they are interpolating
sequences and a different perturbation makes them MZ-arrays. Thus, in a sense,
they behave like the roots of unity in higher dimensions. Since the densities of
the MZ-arrays and the interpolating arrays are well understood, see
\cite[Theorem 1.6]{Marzo07}, then we will get some geometric information on the
Fekete points.

In the next section we will provide the connection of Fekete points and
interpolating and MZ-arrays. On the last section we will draw some
geometric/metric consequences.

\textbf{Acknowledgment} This paper has its origins in a conversation of the
second author with Robert Berman at the Mittag-Leffler Institute over the
possibility of connecting Fekete points and sampling sequences. It is a pleasure
to thank him for sharing his thoughts and to the Institute for the warm
hospitality and great atmosphere. 

\section{Fekete points, MZ-arrays, Interpolating arrays}

\begin{theorem}
Given $\varepsilon>0$ let $L_\varepsilon=[(1+\varepsilon)L]$ and
\[
\Z_\varepsilon (L)=\Z(L_\varepsilon)=\{z_{L_\varepsilon, 1}, \ldots ,
z_{L_\varepsilon, \pi_{L_\varepsilon}}  \},
\]
where $\Z(L)$ is the set of Fekete points of degree $L,$ then $\Z_\varepsilon=\{
\Z_\varepsilon(L) \}_{L\ge 0}$ is an $L^p$-MZ array, for any $1\le p\le \infty$.
\end{theorem}

\begin{proof}
Assume that $\Z$ is a collection of Fekete points. It satisfies a nice
separation property  that is convenient to prove the first inequality of
\eqref{def-MZ}.

\begin{definition}
A triangular array $\Z$ is uniformly separated if there is a positive number
$\varepsilon>0$ such that
\[
d(z_{Lj},z_{Lk})\ge \frac{\varepsilon}{L+1},\text{ if } j\neq k,
\]
for all $L\ge 0$, where $d(z,w)=\arccos\langle z, w\rangle$.
\end{definition}

Reimer, \cite{Reimer90} observed that the Fekete points are uniformly
separated. More precisely,
\[
\frac{\pi}{2L}\le \min_{i\neq j}d(z_{Li},z_{Lj})
\]
just by using Marcel Riesz result for trigonometric polynomials on great
circles.Thus we know that
\[
\min_{i\neq j}d(z_{L_\varepsilon i},z_{L_\varepsilon j})\ge
\frac{\pi}{2L_\varepsilon}\ge \frac{C_\varepsilon}{L+1},
\]
and therefore the array $\Z_\varepsilon$ is uniformly separated. This implies
the following Plancherel-Polya type inequality for any $1\le p<\infty$
\[
\frac{1}{\pi_{L}}\sum_{j=1}^{\pi_{L_\varepsilon}}|Q(z_{L_\varepsilon
j})|^p\lesssim \int_{\S^d}|Q(z)|^p d\sigma(z), \text{ for any } Q\in \Pi_{L},
\]
see \cite[Corollary 4.6]{Marzo07}.

The right hand side inequality in \eqref{def-MZ} is more delicate, we need an
appropriate representation formula for the polynomials in terms of the values at
the points. The most naive approach is to start by the Lagrange interpolation
formula. Let
\[
\ell_{Li}(z)=\frac{\Delta( z_{L1} ,\ldots ,z_{L,i-1},z,z_{L,i+1},\ldots
,z_{L,\pi_L} )}{\Delta( z_{L1} ,\ldots  ,z_{L,\pi_L} )}
\]
then $\sup_{z\in \S^d} | \ell_{Li}(z) |\le 1$ and the Lagrange interpolation
operator defined in $\mathcal{C}(\S^d)$ as
\[
\Lambda_L(f)(z)=\sum_{j=1}^{\pi_L}f(z_{Lj})\ell_{Lj}(z)
\] 
satisfies
\[
\| \Lambda_L(f) \|_{\infty}\le \pi_L \|f\|_\infty.
\]
We want better control of the norms. So we need a slightly bigger set of points
and a weighted representation formula. Let $p$ be a polynomial in one variable
of degree $[L \varepsilon]$ and such that $p(1)=1$. Then given $Q\in \Pi_L$ one
has for a fixed $z\in \S^d$ 
\[
R(w)=Q(w)p(\la z, w \ra)\in \Pi_{L_\varepsilon} 
\]
and therefore we obtain our weighted representation formula:
\[
Q(z)=\sum_{j=1}^{\pi_{L_\varepsilon}}p(\la z, z_{L_\varepsilon ,
j}\ra)Q(z_{L_\varepsilon , j})\ell_{L_\varepsilon , j}(z). 
\]
We define the operator $Q_L$ from $\C^{\pi_{L_\varepsilon}}\to
\Pi_{L_{2\varepsilon}}$ as
\[
Q_L[v](z)=\sum_{j=1}^{\pi_{L_\varepsilon}}v_j p(\la z, z_{L_\varepsilon , j}\ra)
\ell_{L_\varepsilon , j}(z)\qquad \forall v\in \C^{\pi_{L_\varepsilon}}.
\]
We want to prove that 
\begin{equation}\label{intbound}
\int_{\S^d}|Q_L[v](z)|^p\,d\sigma(z) \lesssim \frac 1{\pi_{L_\varepsilon}}
\sum_{j=1}^{\pi_{L_\varepsilon}} |v_j|^p,
\end{equation}
with constants uniform in $L$ which is the righthand sided of \eqref{def-MZ}. We
need to choose the weight $p$ with care. We need a polynomial $p$ that peaks at
one point, has degree $[\varepsilon L]$ and decays fast far away from the
picking point. For this purpose we will use powers of the Jacobi polynomials
which are natural in this context because they are the reproducing kernels in
$\Pi_L$, see \cite{Marzo07}. The Jacobi polynomials $P_{L}^{(\alpha,\beta)}$ of
degree $L$ and index $(\alpha,\beta)$ are the orthogonal polynomials on $[-1,1]$
with respect to the weight function $(1-x)^{\alpha}(1+x)^{\beta}$ with
$\alpha,\beta>-1$. We take the normalization
\[
P_{L}^{(\alpha,\beta)}(1)=\binom{L+\alpha}{L}\simeq L^{\alpha}.
\]

We can use the estimates in \cite[Section 7.34]{Szego75} to obtain, for any
$v\in \S^{d}$
\begin{equation}\label{norma}
\int_{\S^{d}}|P_{L}^{(d/2,d/2-1)}(\la u, v\ra)|^{2}d\sigma(u) \simeq 1,\qquad
\forall L>0.
\end{equation}

We will use as auxiliary polynomial
\[
p(t)=L^{-d}\left( P_{[\varepsilon L/2]}^{(d/2,d/2-1)}(t) \right)^2,
\]
then $p(1)\simeq1$ and by the estimate \eqref{norma}
\[
\int_{\S^d} |p(\la z, z_{L_\varepsilon , j}\ra)| d\sigma (z)\simeq L^{-d}
\int_{\S^d} |P_{[\varepsilon L/2]}^{(d/2,d/2-1)}(\la z, z_{L_\varepsilon ,
j}\ra)|^2 d\sigma (z)\simeq L^{-d}\simeq\pi_{L}^{-1}.
\]

Now as $|\ell_{L_\varepsilon , j}(z)|\le 1$ one has
\begin{align*}
\int_{\S^d} & |Q_L[v]| d\sigma (z) \le \sum_{j=1}^{\pi_{L_\varepsilon}}|v_j|
\int_{\S^d} |p(\la z, z_{L_\varepsilon , j}\ra)| d\sigma (z)\lesssim
\frac{1}{\pi_L}\sum_{j=1}^{\pi_{L_\varepsilon}}|v_j|
\end{align*}
and also for any fixed $z\in \S^d$
\begin{align*}
|Q_L[v](z)| & \le \sup_j |v_j| \sum_{j=1}^{\pi_{L_\varepsilon}} |p(\la z,
z_{L_\varepsilon , j}\ra)| \\ & \le \sup_j |v_j|  \int_{\S^d} \pi_L |p(\la z,
z_{L_\varepsilon , j}\ra)| d\sigma (z) \lesssim \sup_j |v_j|.
\end{align*}
Then the result follows by the Riesz-Thorin interpolation theorem.
\end{proof}

The corresponding result for interpolation reads as follows:

\begin{theorem}
Given $\varepsilon>0$ let $L_{-\varepsilon}=[(1-\varepsilon)L]$ and let
\[
\Z_{-\varepsilon} (L)=\Z(L_{-\varepsilon})=\{z_{L_{-\varepsilon}, 1}, \ldots ,
z_{L_{-\varepsilon}, \pi_{L_{-\varepsilon}}}  \},
\]
where $\Z(L)$ is the set of Fekete points of degree $L,$ then
the array $\Z_{-\varepsilon}=\{ \Z_{-\varepsilon}(L) \}_{L\ge 0}$ is 
$L^p$-interpolating, for any $1\le p\le \infty$.
\end{theorem}

\begin{proof}
Given an array of values $\{ v_{L_{-\varepsilon}j}
\}_{j=1}^{\pi_{L_{-\varepsilon}}}$, we can define the polynomials in $\Pi_L$
\[
R_L[v](z)=\sum_{j=1}^{\pi_{L_{-\varepsilon}}}v_{L_{-\varepsilon}j}p(\la z,
z_{L_{-\varepsilon} , j}\ra)\ell_{L_{-\varepsilon}j}(z),
\]
and $R_L(z_{L_{-\varepsilon}j})=v_{L_{-\varepsilon}j}$. This time the map $R_L$
is from $\C^{\pi_{L_{-\varepsilon}}}\to \Pi_L$ and the $L^p$-estimates on the
norm of $R_L$ follow exactly as the estimates of $Q_L$ in the previous Theorem.
\end{proof}

\section{Geometric properties of the Fekete points}
We will draw some geometric information on the Fekete points. For a given
$z\in\S^d$ and $0<R<1$ we denote by $B(z,R)$ the spherical cap $B(z,R)=\{w\in
\S^d; d(z,w)<R\}$. We will prove that as $L\to\infty$ the number of Fekete
points in $B(z,R)$ gets closer to $\pi_L\tilde\sigma(B(z,R))$  where
$\tilde{\sigma}$ is the normalized Lebesgue measure on $\S^d$, i.e.
$\tilde\sigma=\sigma/\sigma(\S^d)$. We need first information on the density of
M-Z and interpolation arrays.

\begin{definition} For $\Z$ a uniformly separated triangular array in $\S^d$ 
we define the upper and  lower density of the array respectively as
\[
D^{-}(\Z)=\liminf_{\alpha\to\infty}\liminf_{L\to\infty}\frac{\min_{z\in\S^d
}\#(\Z(L) \cap B(z,\alpha/L))/\pi_L}{\tilde\sigma(B(z,\alpha/L))},
\]
\[
D^{+}(\Z)=\limsup_{\alpha\to\infty}\limsup_{L\to\infty}\frac{\min_{z\in\S^d
}\#(\Z(L) \cap B(z,\alpha/L))/\pi_L}{\tilde\sigma(B(z,\alpha/L))}.
\] 
\end{definition}
The main result in \cite[Theorem 1.6]{Marzo07} is (with a slight different
notation)
\begin{theorem} Let $1 \le p \le \infty$. Let $\Z$ be a uniformly separated
array. If $\Z$ is an $L^p$-Marcinkiewicz-Zygmund array then 
$ D^{-}(\Z)\ge 1$. On the other hand if $\Z$ is an $L^p$-interpolating array
then $D^{+}(\Z)\le 1$.
\end{theorem}

Let $\Z(L)$ be the set of Fekete points of degree $L$. We know that for any
$\varepsilon>0$ the array $\Z_\varepsilon=\{ \Z_\varepsilon(L) \}_{L\ge 0}$ is
$L^2$-MZ, so if we unwind the definitions corresponding to the densities, we get
that for any $\varepsilon>0$, there is a big $\alpha=\alpha(\varepsilon)$ such
that for all $L$ and $z\in \S^d$
\begin{equation}\label{d1}
\frac{\frac{1}{\pi_L}\# (\Z(L)\cap
B(z,\frac{\alpha}{L}))}{\tilde{\sigma}(B(z,\frac{\alpha}{L}))}\ge
(1-\varepsilon).
\end{equation}

Similarly since  $\Z_{-\varepsilon}$ is interpolating, whenever $\Z$ is a Fekete
array, from the density condition we get that there is a big
$\alpha=\alpha(\varepsilon)$ such that for all $L$ and $z\in \S^d$
\begin{equation}\label{d2}
\frac{\frac{1}{\pi_L}\# (\Z(L)\cap
B(z,\frac{\alpha}{L}))}{\tilde{\sigma}(B(z,\frac{\alpha}{L}))}\le
(1+\varepsilon).  
\end{equation}
 
\subsection{Vague convergence}

Let us see how the inequalities \eqref{d1} and \eqref{d2} imply that the
normalized counting measure converges vaguely to the Lebesgue measure. Indeed,
defining
\[
\mu_L=\frac{1}{\pi_L}\sum_{j=1}^{\pi_L}\delta_{z_{Lj}}
\]
we have
\begin{align*}
(\mu_L\ast \chi_{B(N,\alpha/L)})(z) & =\int_{\nu \in SO(d+1)}
\chi_{B(N,\alpha/L)}(\nu^{-1}z)\mu_L(\nu N)d\nu \\ & =
\frac{1}{\pi_L}\#(\Z(L)\cap B(z,\alpha/L))
\end{align*}
and
\[
(\sigma \ast \chi_{B(N,\alpha/L)})(z)=\sigma(B(z,\alpha/L))
\]
and finally for any $\varepsilon>0$ there is a big $\alpha$ such that 
\begin{equation}\label{dens-control}
(1-\varepsilon)(\tilde{\sigma} \ast \chi_{B(N,\alpha/L)})(z)\le  (\mu_L\ast
\chi_{B(N,\alpha/L)})(z) \le (1+\varepsilon)(\tilde{\sigma} \ast
\chi_{B(N,\alpha/L)})(z),
\end{equation}
for any $z$ and $L\ge 1$. 

We take an arbitrary spherical cap $B(z,r)$. We want to check that
$\mu_L(B(z,r))\to \tilde{\sigma}(B(z,r))$ as $L\to\infty$.
We fix an $\varepsilon>0$ and we take the convolution of \eqref{dens-control}
with the function $\frac{\chi_{B(N,r)}}{\tilde\sigma(B(N,\alpha/L))}$ with a
very big $L$ and this proves
\[
(1-\varepsilon)(\tilde{\sigma} \ast \chi_{B(N,r-\alpha/L)})(z)\le  (\mu_L\ast
\chi_{B(N,r+\alpha/L)})(z)
\]
and
\[
(\mu_L\ast\chi_{B(N,r-\alpha/L)})(z)\le  (1+\varepsilon)(\tilde{\sigma} \ast
\chi_{B(N,r+\alpha/L)})(z),
\]
for any $z$ and $L$ big. We take limits as $L\to\infty$ and since $\Z$ is
uniformly separated that means that $\mu_L(B(z,r+\alpha/L)\setminus
B(z,r-\alpha/L))\to 0$ uniformly in $z\in\S^d$ as $L\to\infty$. Thus
\[
\lim_{L\to\infty}\mu_L(B(z,r))= \tilde{\sigma}(B(z,r)).
\]
This is for an arbitrary spherical cap. This already implies the vague
convergence of the measures, see \cite{Blumlinger90}, i.e.
\[
\lim_{L\to \infty}
\frac{1}{\pi_L}\sum_{j=1}^{\pi_L}f(z_{Lj})=\frac{1}{\sigma(\S^d)}\int_{\S^d}
f(z)d\sigma (z)
\]
for any $f\in \mathcal{C}(\S^d)$.

From the uniform densities condition on the Fekete points we may obtain other
geometric consequences on the distribution of the Fekete points $\Z(L)$. We
give one example. It is well known that using the bound given by by
L. Fejes T\'oth \cite{Toth49} for the maximum of the minimal spherical distance
between any set of $\pi_{L}=(L+1)^{2}$ points on $\S^{2}$,   there exist at
least two points in $\Z(L)$ at distance $d_{L}$ with
\[
d_{L}\le \arccos \frac{\cot^{2} \omega_{L}-1}{2},\qquad
\omega_{L}=\frac{(L+1)^{2}}{(L+1)^{2}-2}\frac{\pi}{6},
\]
but as
\[
L \arccos \frac{\cot^{2} \omega_{L}-1}{2}\nearrow \kappa= 3.80925\ldots,\quad
\text{ when }L\to \infty \] one has
\[
\min_{i\neq j} d(z_{Li},z_{Lj})\le  \frac{\kappa}{L}.
\]
However the numerical results in \cite{SloWom04} suggest that the right bound
should be 
\[
\limsup_{L\to\infty}L\min_{i\neq j} d(z_{Li},z_{Lj})\le  \pi.
\]

To bound the maximal number, $N$, of disjoints spherical caps in $\S^2$ of
radius $\eta/L$ to be found in a larger spherical cap of radius
$\alpha/L$ we use the following result due to J. Moln\'ar 
\[
N\le
\frac{\pi}{\sqrt{12}}\frac{\sigma(B(z,\frac{\alpha}{L}))}{\sigma(B(z,\frac{
\eta}{2L}))},
\]
see \cite{Molnar52}. Then substituting in the density condition one gets
\[
\eta\le 4 \sqrt{\frac{\pi}{\sqrt{12}}}=\kappa.
\] 
So, we get not only that there exists a couple of points in each generation
$\Z(L)$ at distance smaller than $\kappa/L$ but a more uniform estimate: for any
$\varepsilon>0 $ there is an $\alpha$ such that:
\[ 
\inf_{z_{Li},z_{Lj}\in B(z,\alpha/L)}d(z_{Li},z_{Lj})\le
\frac{\kappa+\varepsilon}{L}
\]
for any spherical cap $B(z,\alpha/L)$.

Up to now we have drawn information from the MZ-arrays and interpolating arrays
to get new information on the Fekete points, but the reverse trend can also be
useful. For instance, since any Fekete array has density one and a small
perturbation makes it a MZ-array (or interpolation), then we obtain the
following corollary
\begin{corollary}
Given any $\epsilon>0$, there are arrays $\Z_{\varepsilon}$ and
$\Z_{-\varepsilon}$ with densities
$D^{+}(\Z_\varepsilon)=D^{-}(\Z_\varepsilon)=1+\varepsilon, D^{+}
(\Z_{-\varepsilon})=D^{-}(\Z_{-\varepsilon})=1-\varepsilon$, such that
$\Z_\varepsilon$ is an $L^p$-MZ array for any $p\in[1,\infty]$ and 
$\Z_{-\varepsilon}$ is an $L^p$-interpolating array for any $p\in[1,\infty]$.
\end{corollary}
Thus the necessary density conditions that $MZ$-arrays and interpolating arrays
satisfy are sharp.

\end{document}